\documentclass{article}

\usepackage{amsmath}
\usepackage{amssymb}
\usepackage{latexsym}
\usepackage{amsxtra}
\usepackage{amscd}
\usepackage{theorem}
\usepackage{xypic}

\theoremstyle{plain}

{\theorembodyfont{\slshape}

        \newtheorem{thm}{Theorem}[section]
        \newtheorem{cor}[thm]{Corollary}
        \newtheorem{lem}[thm]{Lemma}
        \newtheorem{prop}[thm]{Proposition}

}

{\theorembodyfont{\rmfamily}

        \newtheorem{defn}[thm]{Definition}
        
        \newtheorem{rem}[thm]{Remark}
        \newtheorem{exa}[thm]{Example}
        
        \newtheorem{notation}[thm]{Notation}

}

\renewcommand{\em}{\sl}

\newcommand{\proof}{{\bf Proof:\ }}
\newcommand{\Endproof}{\hspace*{\fill} $\Box$ \vspace{1ex} \noindent }

\makeatletter
\renewcommand{\subsubsection}{\@startsection{subsubsection}{3}%
        {\z@}{-3.25ex plus -1ex minus-.2ex}{-1em}{\bf}}
\makeatother

\renewcommand{\arraystretch}{1.5}

\newcommand{\ZZ}{\mathbb{Z}}
\newcommand{\CC}{\mathbb{C}}
\newcommand{\RR}{\mathbb{R}}
\newcommand{\QQ}{\mathbb{Q}}

\newcommand{\F}{\mathcal{F}}

\newcommand{\U}{\mathcal{U}}

\newcommand{\X}{\mathcal{X}}
\newcommand{\Y}{\mathcal{Y}}

\newcommand{\Xs}{\mathsf{X}}
\newcommand{\Ys}{\mathsf{Y}}
\newcommand{\As}{\mathsf{A}}
\newcommand{\Bs}{\mathsf{B}}
\newcommand{\Cs}{\mathsf{C}}
\newcommand{\Us}{\mathsf{U}}

\newcommand{\Ds}{\mathsf{D}}

\newcommand{\OO}{\mathcal{O}}
\newcommand{\HH}{\mathcal{H}}

\newcommand{\ac}{^{\rm ac}}
\newcommand{\nr}{^{\rm nr}}

\newcommand{\ad}{^{\rm ad}}

\newcommand{\univ}{^{\rm uv}}

\newcommand{\para}{_{\rm p}}

\newcommand{\Fh}{\hat{F}}
\newcommand{\Kh}{\hat{K}}
\newcommand{\OOh}{\hat{\OO}}

\newcommand{\QQb}{\bar{\QQ}}

\newcommand{\Hom}{{\rm Hom}}

\newcommand{\Gal}{{\rm Gal}}

\newcommand{\GL}{{\rm GL}}
\newcommand{\SL}{{\rm SL}}

\newcommand{\p}{\mathfrak{p}}
\newcommand{\Spec}{{\rm Spec\,}}
\newcommand{\Spf}{{\rm Spf}}
\newcommand{\Spm}{{\rm Spm}}
\newcommand{\Frac}{{\rm Frac}}

\newcommand{\sw}{{\rm sw}}

\newcommand{\Ind}{{\rm Ind}}
\newcommand{\cInd}{\text{\rm c-Ind}}
\newcommand{\val}{{\rm val}}
\newcommand{\Log}{{\rm Log\,}}

\newcommand{\inj}{\hookrightarrow}
\newcommand{\To}{\;\longrightarrow\;}

\newcommand{\gen}[1]{\mathopen\langle#1\mathclose\rangle}
\newcommand{\abs}[1]{\lvert#1\rvert}
\newcommand{\Abs}[1]{\Bigl\lvert#1\Bigr\rvert}

\makeatletter
\renewcommand{\subsection}{\@startsection{subsection}{2}%
        {\z@}{-3.25ex plus -1ex minus-.2ex}{-1em}{\bf}}
\makeatother

\begin{document}

\title{Swan conductors on the boundary of Lubin-Tate spaces}

\author{Stefan Wewers}

\date{}

\maketitle

\begin{abstract}
  Lubin-Tate spaces of dimension one are finite \'etale covers of the
  non-archimedian open unit disk. We compute certain invariants which
  measure the ramification of this cover over the boundary of the disk. 
\end{abstract}

\section{Introduction}

\subsection{}

Let $F$ be a local field, i.e.\ a finite extension of $\QQ_p$ or the
field of Laurent series over a finite field. Let $\OO$ denote the ring
of integers of $F$. Choose a uniformizer $\varpi$ of $\OO$. The
residue field $k=\OO/(\varpi)$ is a finite field with, say, $q$
elements. We let $\Fh\nr$ denote the completion of the maximal
unramified extension of $F$ and $\OOh\nr\subset\Fh\nr$ its ring of
integers. Fix an integer $h\geq 1$.

A construction due to Lubin-Tate and Drinfeld attaches to $F$ and $h$
a certain inverse system $\cdots\to\X(\varpi^2)\to\X(\varpi)\to\X(1)$
of formal schemes $\X(\varpi^n)=\Spf A_n$, called the {\em Lubin-Tate
tower}. In this paper we only consider the case $h=2$, and we look at
the tower of rigid analytic spaces
\[
      \Xs(\varpi^n):=\X(\varpi^n)\otimes_{\OOh\nr}\Fh\nr,
      \;\;\;n=1,2,\ldots
\]
associated to $(\X(\varpi^n))_n$. At the lowest level, $\Xs(1)$ is
isomorphic to the open unit disk over the non-archimedian field
$\Fh\nr$. The maps $\Xs(\varpi^n)\to\Xs(1)$ are finite \'etale Galois
covers with Galois group $\GL_2(\OO/(\varpi^n))$. Therefore, the rigid
space $\Xs(\varpi^n)$ is a smooth analytic curve.

The present paper is concerned with computing certain invariants which
measure the ramification of the \'etale cover $\Xs(\varpi^n)\to\Xs(1)$
over the `boundary' of the open disk $\Xs(1)$. These computations are
a key ingredient for the results of \cite{LT}, which describe the
{\em stable reduction} of $\Xs(\varpi^n)$.

Our notation $\Xs(\varpi^n)$ (which is not standard) suggests an
analogy with modular curves. And indeed, for $F=\QQ_p$ the connected
components of $\Xs(p^n)\otimes_{\QQ_p\nr}\CC_p$ are isomorphic to
certain open analytic subspaces of the classical modular curve
$X(p^nm)$ associated to the full congruence subgroup of $\SL_2(\ZZ)$
modulo $p^nm$, for $m\geq 3$ and prime to $p$. Each of these analytic
subspaces is the formal fiber of a singularity of the $p$-adic
integral model of $X(p^nm)$ studied by Katz and Mazur
\cite{KatzMazur}. Combining the results of \cite{LT} with
the results of Katz and Mazur yields, in the special case $F=\QQ_p$, a
description of the stable reduction of $X(p^nm)$ at the prime $p$.
There is also an interesting connection with the local Langlands
correspondence. A theorem of Carayol (generalized by Harris and Taylor
to the case $h>2$) asserts that the \'etale cohomology of the tower
$\varprojlim \Xs(\varpi^n)$ realizes the local Langlands
correspondence for the group $\GL_2(F)$. Using the description of the
stable reduction of $\Xs(\varpi^n)$ in \cite{LT}, it is possible to
give a new and more direct proof of this theorem (work in progress).

\subsection{}

To explain the result of this paper and their relevance for the
results of \cite{LT}, we abstract a bit from the special situation of
Lubin-Tate spaces. Let $K$ be a complete discrete non-archimedian
field. Let $R$ denote the ring of integers of $K$. Let $\Ys$
denote the open unit disk over $K$ and let $f:\Xs\to\Ys$ be a
finite \'etale Galois cover with Galois group $G$. In what follows, we
shall always allow the field $K$ to be replaced by a suitable finite
extension. (We do not replace $K$ by its algebraic closure because 
we want the valuation to be discrete. See \cite{open}.)

By the semistable reduction theorem, the rigid space $\Xs$ has a
minimal semistable formal model $\X$ over $R$. Let $\Y:=\X/G$ be the
quotient by the $G$-action; this is a semistable formal model of
$\Ys$. Such a semistable model of the disk $\Ys$ is easy to describe:
it is determined by a finite collection of closed affinoid disks
$\Ds_i\subset\Ys$. Namely, each closed disk $\Ds_i$ determines a
blowup of the standard formal model $\Spf R[[T]]$ of $\Ys$, and
performing all these blowups simultaneously yields the semistable
model $\Y$. We say that the disks $\Ds_i$ are {\em relevant} for the
stable reduction of $\Xs$. In some sense, the main problem in finding
the stable reduction of $\Xs$ is to find the relevant disks
$\Ds_i$. See \cite{open} for more details.

Let $\tau$ be an irreducible representation of the group $G$. We
assume that $\tau$ is defined over a finite extension of $\QQ_\ell$,
where $\ell$ is an auxiliary prime dividing neither the order of $G$
nor the characteristic of the field $k$. The representation $\tau$ and
the Galois cover $f:\Xs\to\Ys$ determine a lisse $\QQb_\ell$-sheaf
$\F$ on the \'etale topology of $\Ys$.

Let $x\in\Ys$ be a closed point. For a positive rational number
$s\in\QQ_{>0}$ we let $\Ds(x,s)\subset\Ys$ denote the closed affinoid
disk with center $x$ and radius $r:=\abs{\varpi}^s$. Following Huber
\cite{Huber01} and Ramero \cite{Ramero}, we define two numbers,
\[
   \sw_{\F}(s)\in\ZZ_{\geq 0}, \qquad \delta_{\F}(s)\in\QQ_{\geq 0},
\]
which we call the {\em Swan conductor} and the {\em discriminant
conductor}. These numbers measure, in some sense, the ramification of
the sheaf $\F$ over the disk $\Ds(x,s)$. For instance, the
discriminant conductor $\delta_{\F}(s)$ is essentially the classical
Swan conductor of the Galois representation associated to $\F$ and the
discrete valuation of the field $\Frac(R[[T]])$ corresponding to the
maximum norm on $\Ds(x,s)$.

The functions $s\mapsto\sw_{\F}(s)$ and $s\mapsto\delta_{\F}(s)$
extend to functions on the interval $[0,\infty)\subset\RR$ with the
following properties: (a) $\sw_{\F}$ is right continuous, locally
constant and decreasing, (b) $\delta_{\F}$ is continuous and piecewise
linear, and (c) $\sw_{\F}$ is equal to minus the right derivative of
$\delta_{\F}$,
\[
      \frac{\partial}{\partial s}\delta_{\F}(s+) = -\sw_{\F}(s).
\]
It follows that $\delta_{\F}$ is convex and decreasing and that both
$\sw_{\F}$ and $\delta_{\F}$ are $\equiv 0$ for $s\gg 0$. Furthermore,
there are are finite number of (rational) critical values
$s_1,\ldots,s_n$ where $\sw_{\F}$ is discontinuous and where
$\delta_{\F}$ is not smooth. We call $s_1,\ldots,s_n$ the {\em breaks}
of $\delta_{\F}$.

The study of the function $\delta_{\F}$ goes back to L\"utkebohmert's
paper \cite{Luetkebohmert93} on the $p$-adic Riemann's Existence
Theorem.  (In {\em loc.cit} one studies finite \'etale covers
$f:\Xs\to\Ys^*:=\Ys-\{0\}$ of the {\em punctured} disk, and the
representation $\tau$ is the regular representation of $G$. In this
case the integers $\sw_{\F}(s)$ may become negative, and $\delta_F(s)$
may never be zero.) To prove the properties of the function
$\delta_{\F}$ mentioned above one uses the semistable reduction
theorem. As an immediate consequence of this proof, one obtains the
following characterization of the breaks $s_1,\ldots,s_n$: the disk
$\Ds(x,s_i)$ are relevant for the stable reduction of the cover
$f:\Xs\to\Ys$.  This is the reason why we are interested in studying
the function $\delta_{\F}$.

\subsection{}

To apply the previous discussion to the Lubin-Tate tower, we set
$\Ys:=\Xs(1)$, $\Xs:=\Xs(\varpi^n)$ and
$G=G_n:=\GL_2(\OO/\varpi^n)$. We also have to choose a point
$x\in\Xs(1)$ (the center for the disks we look at) and a
representation $\tau$ of $G_n$ (which gives rise to the sheaf $\F$). It
turns out that, in order to describe the semistable reduction of
$\Xs(\varpi^n)$, it suffices to choose the following pairs $(x,\tau)$:
\begin{itemize}
\item
  The point $x\in\Xs(1)$ is a {\em canonical point} corresponding to a
  separable quadratic extension $E/F$ (in the case $F=\QQ_p$ such a
  point gives rise to a CM-point (of a certain kind) on the
  corresponding modular curve).
\item
  The restriction of the representation $\tau$ to the compact group
  $K:=\GL_2(\OO)$ is a {\em type} for an irreducible supercuspidal
  representation of $\GL_2(F)$. Moreover, the quadratic extension
  associated to $\tau$ by the classification of such types (see e.g.\
  \cite{Kutzko78} or \cite{HenniartAppendix}) is the extension $E/F$.
\end{itemize}
For each pair $(x,\tau)$ as above, one can explicitly compute the
function $\delta_{\F}$ and, in particular, determine its breaks
$s_1,\ldots,s_r$. As we have remarked above, the disks $\Ds(x,s_i)$ are
all relevant for the stable reduction of $\Xs(\varpi)$. It is not true
that every relevant disk occurs in this way, but nevertheless, knowing
these ones is sufficient to study the stable reduction of
$\Xs(\varpi^n)$. See \cite{LT}. 

We hope that these explanations sufficiently motivate the main result
of this paper. It computes the value of functions $\sw_{\F}$ and
$\delta_{\F}$ at the point $s=0$ in all the cases we are interested
in.

\begin{thm}
  Let $\tau$ be the $K$-type of an irreducible supercuspidal
  representation of $\GL_2(F)$. Assume that $\tau$ is minimal of level
  $n$, and consider $\tau$ as a representation of
  $G_n=\GL_2(\OO/\varpi^n)$. (See \S \ref{supercuspidal2} for the
  terminology.) Let $\F$ be the sheaf on $\Xs(1)$ corresponding to
  $\tau$. Let $E/F$ be the quadratic extension associated to $\tau$ by
  the construction in \cite{Kutzko78}. Let $x\in\Xs(1)$ be a canonical
  point, corresponding to the extension $E/F$ (actually, {\em any}
  point $x\in\Xs(1)$ would do). If the extension $E/F$ is unramified
  then we have
  \[
      \sw_{\F}(0) = -(q+1)q^{n-1}, \quad
      \delta_{\F}(0) = (nq-n+1)q^{n-1}.
  \]
  If $E/F$ is ramified then
  \[
      \sw_{\F}(0) = -(q+1)q^{n-2}, \quad
      \delta_{\F}(0) = (nq-q-n)q^{n-2}.
  \]
\end{thm}

Note that, a priori, the values of $\sw_{\F}$ and $\delta_{\F}$ at
$s=0$ are not defined in the same way as for rational values $s>0$,
simply because the closed unit disk $D(x,0)$ is not contained in the
open unit disk $\Xs(1)$. There are essentially two ways to get hold of
the value at $s=0$. Firstly, one can try to extend the cover
$f:\Xs(\varpi^n)\to\Xs(1)$ to the closed unit disk. This is possible
because the cover $f$ essentially occurs inside an algebraic cover
between certain Shimura or Drinfeld modular curves (classical modular
curves for $F=\QQ_p$). However, this is not very helpful to actually
compute $\sw_{\F}(s)$ and $\delta_{\F}(s)$ for $s=0$. So we follow the
second possibility, i.e.\ we extend the definition of $\sw_{\F}$ and
$\delta_{\F}$ in a way that allows us to compute these values for
$s=0$ directly in terms of the rank-two valuation of $\Frac(R[[T]])$
corresponding to the 'boundary' of the open unit disk.


\section{The Swan and the discriminant conductor}

We introduce Swan and discriminant conductors for \'etale Galois
covers of so-called {\em open analytic curves}, extending results
of Huber \cite{Huber01} and Ramero \cite{Ramero}, where this is done
for smooth affinoid curves.

\subsection{Recall: open analytic curves}

We fix a field $K_0$ which is complete with respect to a discrete
non-archimedian valuation $\abs{\,\cdot\,}$, and whose residue field
$k$ is algebraically closed and of positive characteristic $p>0$. We
choose an algebraic closure $K\ac$ of $K_0$ and extend the valuation
$\abs{\,\cdot\,}$ to $K\ac_0$. We let $\Gamma$ denote the value group
of the valuation $\abs{\,\cdot\,}$ on $K\ac$. We assume that $\Gamma$
is a subgroup of $\RR_{>0}$.

We recall the following definitions from \cite{open}.

\begin{defn} \label{opendef}
  An {\em open analytic curve} is given by a pair $(K,\Xs)$, where
  $K\subset K\ac_0$ is a finite extension of $K_0$ and $\Xs$ is a
  rigid analytic space over $K$. We demand that $\Xs$ is isomorphic to
  $C-D$, where $C$ is the analytification of a smooth projective curve
  over $K$ and $D\subset C$ is an affinoid subdomain intersecting
  every connected component of $C$.

  A morphism between two open analytic curves $(K_1,\Xs_1)$ and $(K_2,\Xs_2)$
  is an element of the direct limit 
  \[
        \Hom(\Xs_1,\Xs_2):= \varinjlim_{K_3} \Hom(\Xs_1\otimes K_3,\Xs_2\otimes K_3),
  \]
  where $K_3\subset K\ac_0$ ranges over all common finite extensions
  of $K_1$ and $K_2$.
\end{defn}

We shall follow the same convention as in \cite{open}, regarding the
field $K$. That is, we simply write $\Xs$ to denote an open analytic
curve. The field $K$ is only mentioned if needed, and it is always
assumed to be `sufficiently large'. For instance, if we say that $\Xs$
is connected then this means that $\Xs\otimes_K K'$ is connected, for every
finite extension $K'/K$.

\begin{defn} \label{opendef2} An {\em underlying affinoid} is an
  affinoid subdomain $\Us\subset\Xs$ such that $\Xs-\Us$ is the disjoint
  union of open annuli none of which is contained in an affinoid subdomain
  of $\Xs$. An {\em end} of $\Xs$ is an element of the inverse limit
  of the set of connected components of $\Xs-\Us$, where $\Us$ ranges over
  all underlying affinoids. The set of all ends is denoted by
  $\partial\Xs$.
\end{defn}

\subsection{The rank two valuation associated to an end}

Let $\Xs$ be an open analytic curve with field of definition $K$. We
let $R$ denote the valuation ring of $K$. Let $x\in\partial \Xs$ be an
end of $\Xs$.  We will associate to the pair $(\Xs,x)$ a certain field
$K_x$ equipped with a rank two valuation $\abs{\,\cdot\,}_x$. The
construction should depend only on the connected component of $\Xs$ on
which $x$ lies, so we may from the start assume that $\Xs$ is
connected.

Recall from \cite{open} that $\Xs$ has a certain canonical formal
model $\X$ called the {\em minimal model}, as follows. Let
$\OO_{\Xs}^\circ$ denote the ring of power bounded analytic functions
on $\Xs$.  According to \cite{Bosch77}, $\OO_{\Xs}^\circ$ is a normal
and complete local ring, isomorphic to the completion of an
$R$-algebra of finite type. We set $\X:=\Spf\, \OO_{\Xs}^\circ$; then
the generic fiber of $\X$ can be canonically identified with $\Xs$. By
choosing the field $K$ sufficiently large we may assume that the
scheme $\X_s:=\Spec(\OO_{\Xs}^\circ\otimes k)$ is reduced. Under this
assumption, the construction of $\X$ is stable under further extension
of $K$, i.e.\ if $K'/K$ is a finite extension with valuation ring $R'$
then $\X'=\X\otimes_R R'$ is the minimal model of $X\otimes_K K'$.

There is a natural bijection between the ends of $\Xs$ and the generic
points of the scheme $\X_s$. We let $\eta\in\X_s$ denote the generic
point corresponding to $x\in\partial\Xs$. To $\eta$ corresponds a
discrete valuation of the fraction field of $\OO_{\Xs}^\circ$ which
extends the valuation $\abs{\,\cdot\,}$ on $K$. We extend this
valuation to the field $\Frac(\OO_{\Xs}^\circ\otimes_R K\ac)$ in the
obvious way and denote it by $\abs{\,\cdot\,}_\eta$. The residue field
$k(\eta)$ of $\abs{\,\cdot\,}_\eta$ is a field of Laurent series,
$k(\eta)\cong k((t))$. We define the rank two valuation
$\abs{\,\cdot\,}_x$ on $\Frac(\OO_{\Xs}^\circ\otimes_R K\ac)$ as the
composition of $\abs{\,\cdot\,}_\eta$ with the canonical valuation on
the residue field $k(\eta)$ (see e.g.\ \cite{ZariskiSamuel}). We
define the field
\[
      K_x := \Frac(\OO_{\Xs}^\circ\otimes_R K\ac)\sphat
\]
as the completion of $\Frac(\OO_{\Xs}^\circ\otimes_KK\ac)$ with respect to
$\abs{\,\cdot\,}_x$. One easily checks that the field $K_x$ is henselian.

\begin{notation}
Let $\Gamma_x$ denote the valuation group of $\abs{\,\cdot\,}_x$. Then
\begin{equation} \label{Gammaeq}
    \Gamma_x = \Gamma \times \Lambda_x, \quad \Lambda_x=\gen{\gamma_x},
\end{equation}
where $\Lambda_x$ is an ordered cyclic group generated by an element
$\gamma_x<1$. The ordering on $\Gamma_x$ is lexicographic, i.e.\ for
$r,s\in\Gamma$ and $i,j\in\ZZ$ we have $(r,\gamma_x^i)<(s,\gamma_x^j)$
if either $r<s$, or if $r=s$ and $i>j$. Let $\abs{\,\cdot\,}_x^\flat$
denote the rank one valuation on $K_x$ which is the composition of
$\abs{\,\cdot\,}_x$ with the first projection $(r,\gamma_x^i)\mapsto
r$ (this was denoted $\abs{\,\cdot\,}_\eta$ before). Let
$\#:\Gamma_x\to\ZZ$ be the group homomorphism $(r,\gamma_x^i)\mapsto
i$ and write $\#_x$ for the composition of $\abs{\,\cdot\,}_x$ with
$\#$.
\end{notation}

The splitting \eqref{Gammaeq} of $\Gamma_x$ is canonical. To see this,
note that we have identified $\Gamma$ with the maximal divisible
subgroup of $\Gamma_x$. Furthermore, the generator $\gamma_x$ of
$\Lambda_x$ is the maximal element of $\Gamma_x$ which is strictly
smaller than $1$.  An element $u\in K_x$ such that
$\abs{u}_x=\gamma_x$ is called a {\em parameter} for the end $x$.

\begin{exa} \label{annulusexa}
  Let $\epsilon\in \Gamma$, $\epsilon<1$. Choose a finite extension
  $K/K_0$ such that $\epsilon\in\abs{K^\times}$. We regard the
  standard open annulus
  \[
       \Xs := \Cs(\epsilon,1)=\{\, u \,\mid\, \epsilon<\abs{u}< 1\,\}
  \]
  as an open analytic curve with field of definition $K$. Clearly,
  $\Xs$ has two ends. We let $x\in\partial\Xs$ be the end
  corresponding to the family of open annuli
  $\Cs(\epsilon',1)\subset\Xs$ for $\epsilon<\epsilon'<1$.
  Choose an element $\pi\in R$ with $\abs{\pi}=\epsilon$. Then 
  \[
       \OO_{\Xs}^\circ = R[[\,u,v \,\mid\, uv=\pi\,]].
  \]
  An element of $\OO_{\Xs}^\circ$ can be written as a Laurent series $\sum_i
  c_iu^i$ with $c_i\in K$ such that $\abs{c_i}\leq 1$ for $i\geq 0$
  and $\abs{c_i}\leq \epsilon^{-i}$ for $i<0$. It follows that every
  element of $\Frac(\OO_{\Xs}^\circ\otimes_R K\ac)$ can also be written as a
  Laurent series whose coefficients all lie in a finite extension of
  $K$ and satisfy certain growth conditions. With this notation, we have
  \[
       \abs{\sum_ic_i u^i}_x = \max_{i\in\ZZ}\, (|c_i|,\gamma_x^i).
  \]
  In other words, $\abs{\sum_i c_iu^i}_x^\flat=\max_i\abs{c_i}$, and
  $\#(\sum_i c_iu^i)$ is the first index $i$ where $\abs{c_i}$ takes
  its maximal value.
\end{exa}

\subsection{Functoriality}

Let $f:\Xs\to\Ys$ be an analytic map between open analytic curves.
Given $x\in\partial\Xs$ and $y\in\partial\Ys$, the notation $f(x)=y$
means the following. For every sufficiently small open annulus
$\As\subset\Xs$ representing the end $x$ the restriction of $f$ to
$\As$ is finite onto its image, and $f(\As)\subset\Ys$ is an open
annulus representing $y$. We shall also write $f:(\Xs,x)\to(\Ys,y)$
for a map $f$ satisfying the above condition.

\begin{prop} \label{functorprop}
  The map $f:(\Xs,x)\to(\Ys,y)$ induces a finite extension of valued fields
  \[
         f_x^*:K_y \inj K_x.
  \]
  For all $u\in K_y^\times$ we have
  \[
        \abs{f_x^*u}_x^\flat = \abs{u}_y^\flat
  \]
  and 
  \[
        \#_x(f_x^*u) = [K_y:K_x]\cdot\#_y(u).
  \]
  The index $[K_y:K_x]$ is equal to the degree of $f|_{\As}:\As\to
  f(\As)$, where $\As\subset\Xs$ is a sufficiently small annulus
  representing the end $x$.
\end{prop}
 
\proof 
We prove the proposition first in two special cases.

\vspace{1ex}\noindent
{\bf Case 1:} $\Xs$ is an open annulus, $f:\Xs\inj\Ys$ is an
open immersion and $\Xs\subset\Ys$ represents the end $y$.

By \cite{open}, \S 1, we can represent $\Ys$ as the formal fiber
$]z[_Y$, where $Y=\Spec A$ is an affine normal curve over $R$ with
reduced special fiber $Y_s$ and $z\in Y_s$ is a closed point.  The end
$y$ corresponds to a branch $\eta$ of $Y_s$ through $z$.  Here we
consider the branch $\eta$ as a generic point of 
$\Spec(\OOh_{Y,z}\otimes k)$. Let $Z\subset Y_s$ denote the
irreducible component whose generic point is the image of $\eta$ on
$Y_s$.

Furthermore, we can identify the open subspace $\Xs\subset\Ys$ with
the formal fiber $]w[_X$, where $g:X\to Y$ is an admissible blowup
with center $z$ and $w\in X_s$ is a closed point. The assumption that
$\Xs$ is an open annulus representing the end $y$ implies that $w$ is
a formal double point, point of intersection of the strict transform
of $Z$ with the exceptional divisor of $g$. Let
$\nu\in\Spec(\OOh_{X,w}\otimes k)$ denote the branch of $X_s$ through
$w$ corresponding to the strict transform of $Z$. To prove the
proposition in Case 1 it suffices to show that the map
\[
      g_x^*:\OO_{\Ys}^\circ = \OOh_{Y,z} \To \OO_{\Xs}^\circ = \OOh_{X,w}
\]
induced by $g$ is compatible with the valuations $\abs{\,\cdot\,}_y$
and $\abs{\,\cdot\,}_x$ (i.e.\ $\abs{g_x^*u}_x=\abs{u}_y$) and induces
an isomorphism $K_y\cong K_x$.

The valuations $\abs{\,\cdot\,}_y^\flat$ and
$\abs{\,\cdot\,}_x^\flat$, restricted to $\Frac(\OOh_{Y,z})$ and
$\Frac(\OOh_{X,w})$, are discrete and correspond to the codimension one
points $\eta\in\Spec(\OOh_{Y,z})$ and $\nu\in\Spec(\OOh_{X,w})$. Since
$g(\nu)=\eta$ and $X_s$ and $Y_s$ are reduced we have
$\abs{g_x^*u}_x^\flat=\abs{u}_y^\flat$ for all $u\in\OOh_{Y,z}$. The
functions $\#_y$ and $\#_x$, restricted to $\Frac(\OOh_{Y,z})^\times$
and $\Frac(\OOh_{X,w})^\times$, are induced from the natural discrete
valuations on the residue fields of $\eta$ and $\nu$. It is also clear
that $g_x^*$ induces an isomorphism between these two residue fields.
Therefore, $\abs{g_x^*u}_x=\abs{u}_y$ holds for all $u\in\OOh_{Y,z}$.
It remains to show that $g_x^*$ induces an isomorphism $K_y\cong K_x$.
Actually, it suffices to show that $\OOh_{X,w}\subset g_x^*K_y$. Let
$u\in\OO_{Y,z}$ be a parameter for $y$, i.e.\ an element with
$\abs{u}_y=\gamma_y$. By the definition of the valuation
$\abs{\,\cdot\,}_y$, $u$ lies in the maximal ideal of $\OOh_{Y,z}$, so
we have $R[[u]]\subset \OOh_{Y,z}$. Identify the ring $R[[u]]$ with
its image in $\OOh_{X,w}$ via $g_x^*$. Then $u$ is also a parameter
for $\abs{\,\cdot\,}_x$. Therefore,
\[
      \OOh_{X,w} = R[[\,u,v\,\mid\, uv=\pi\,]],
\]
where $\pi\in R$ is a suitable element with $\epsilon:=\abs{\pi}<0$. Consider
$v=\pi/u$ as an element in $K_y$. Using
\[
      \abs{v^i}_y = (\epsilon^i,\gamma_y^{-i}) \To 0
\]
for $i\to\infty$, we conclude that $\OOh_{X,w}$ is contained in $K_y$,
as desired. This completes the proof of the proposition in Case 1.

\vspace{1ex}\noindent
{\bf Case 2:} $\Xs$ and $\Ys$ are open annuli and $f$ is a finite morphism.

In this case the induced map $f^*:\OO_{\Ys}^\circ\to\OO_{\Xs}^\circ$
is a finite extension of normal local rings. With similar arguments as
in Case 1, one shows that it induces an extension of valued fields 
\[
      f_x^*:(\Frac(\OO_{\Ys}^\circ\otimes_R K\ac),\abs{\,\cdot\,}_y) \to 
        (\Frac(\OO_{\Xs}^\circ\otimes_R K\ac),\abs{\,\cdot\,}_x)
\]
such that $\abs{g_x^*u}_x^\flat=\abs{u}_y^\flat$. It is also easy to
see that the induced extension of residue fields of
$\abs{\,\cdot\,}_y$ and $\abs{\,\cdot\,}_x$ (which are discretely
valued fields) is totally ramified of degree $\deg(f)$. It follows
that the induced extension $K_y\to K_x$ on the completions is of the
same degree and has all the claimed properties. This proves the
proposition in Case 2.

The general case of the proposition follows easily from Case 1 and Case 2.
\Endproof

\subsection{Higher ramification groups}

Let $\Xs$ be an open analytic curve and $x\in\partial\Xs$ an end.  Let
$G$ be a finite group acting faithfully on $\Xs$ and fixing $x$. By
Proposition \ref{functorprop} the action of $G$ on $\Xs$ induces an action of $G$
on the valued field $(K_x,\abs{\,\cdot\,}_x)$. It is easy to see that
$G$ acts faithfully on $K_x$. Following Huber \cite{Huber01} we define
a filtration $(G_h)_{h\in\Gamma_x}$ of higher ramification groups on
$G$.

For any $h\in\Gamma_x$, we set
\[
     G_h \,:=\, \{\,\sigma\in G \,\mid\, h_x(\sigma)\leq h\,\},
\]
where
\[
      h_x(\sigma) \,:=\, \min\,\{\,h\,\mid\,|u-\sigma(u)|_x\leq h\cdot|u|_x\,
                     \forall u\in K_x\,\}.
\] 
It is shown in \cite{Huber01}, Lemma 2.1, that
\begin{equation} \label{hsigmaeq}
      h_x(\sigma) \,=\, \Big|\frac{t-\sigma(t)}{t}\Big|_x,
\end{equation}
where $t\in K_x$ is any {\em parameter} at $x$, i.e.\ an element of $K_x$
with $|t|_x=\gamma_x$. It is easy to see that the group
\[
       P := \cup_{h<1} G_h
\]
is the maximal $p$-subgroup and that $G/P$ is cyclic of order prime to $p$. Let
\[
        h_1>\cdots>h_l
\]
be the elements of the set $\{h(\sigma)\mid \sigma\in G,\sigma\neq 1\}$ which
are $\neq 1$. Set $h_0:=1$. By definition we have
\[
     G=G_{h_0}\supsetneq G_{h_1}\supsetneq \ldots\supsetneq G_{h_l} 
         \supsetneq \{1\}.
\]
The elements $h_i\in\Gamma_x$ for $i\geq 1$ are called the {\em
  jumps} in the filtration $(G_h)_h$.

There is also an upper numbering. Let $\Ys:=\Xs/G$ and let
$y\in\partial\Ys$ be the image of $x$. By Proposition
\ref{functorprop} the extension $K_x/K_y$ is of degree $|G|$. Since
$G$ acts faithfully on $K_x$ and fixes $K_y$, the extension $K_x/K_y$
is actually Galois, with Galois group $G$. Let
$\varphi_{K_x/K_y}:\Gamma_x\otimes\QQ\to\Gamma_y\otimes\QQ$ be the
function defined in \cite{Huber01}, \S 2. For
$\gamma\in\Gamma_y\otimes\QQ$ set $G^\gamma:=G_h$, where
$h:=\varphi_{K_x/K_y}^{-1}(\gamma)$. The elements
\[
   \gamma_i := \varphi_{K_x/K_y}(h_i) \in\Gamma_y\otimes\QQ, \quad i=1,\ldots,r,
\]
are called the {\em jumps in the upper numbering}. Explicitly, we have 
\begin{equation} \label{gammaieq}
  \gamma_i^\flat = 
    \prod_{j=1}^l \Bigl(\frac{h_j^\flat}{h_{j-1}^\flat}\Bigr)^{\abs{G_{h_j}}},\quad
   \#\gamma_i = \sum_{j=1}^i\frac{\#h_j-\#h_{j-1}}{(G:G_{h_j})},
\end{equation}
with $h_j=(h_j^\flat,\#h_j)$ and $\gamma_j=(\gamma_j^\flat,\#\gamma_j)$.

\subsection{The Swan conductor and the discriminant conductor}

We fix an open analytic curve $\Ys$ and a prime number $\ell$
different from $p={\rm char}(k)$. Let $\F$ be a lisse sheaf of
$\QQb_\ell$-vectorspaces on $\Ys$. We say that $\F$ is {\em
  admissible} if there exists a finite group $G$ of order prime to
$\ell$, an \'etale $G$-torsor $f:\Xs\to\Ys$ and a representation
$\tau$ of $G$ on a finite-dimensional $\QQb_\ell$-vector
space $W$ such that 
\[
   \F \cong (f_*\QQb_\ell)[\tau]:=
     \underline{\Hom}_G(\underline{W},f_*\QQb_\ell).
\]
See \cite{open} for more details and arguments why it makes sense to
work with $\QQb_\ell$-coefficients in this situation. 

Let $\F$ be an admissible sheaf on $\Ys$ and $y\in\partial Y$ an end.
We will attach to $\F$ and $y$ two invariants, the {\em Swan
  conductor} $\sw_y(\F)\in\ZZ$ and the {\em discriminant conductor}
$\delta_y(\F)\in\RR$.

Choose an \'etale $G$-torsor $f:\Xs\to\Ys$ and a representation
$\tau:G\to\GL(V)$ such that $\F\cong(f_*\QQb_\ell)[\tau]$. Choose also
an end $x\in\partial\Xs$ with $f(x)=y$. Let $G_x\subset G$ denote the
stabilizer of $x$. In the previous subsection we defined a function
$\sigma\mapsto h_x(\sigma)$ on $G_x$ with values in $\Gamma_x$. We now define
\[
   \sw_x(\sigma) := -\#h_x(\sigma)\;\;\;\text{if $\sigma\neq 1$}, \quad
   \sw_x(1) := \sum_{\sigma\neq 1} \#h_x(\sigma)
\]
and
\[
    \sw_y:= \Ind_{G_x}^{G} \sw_x.
\]
This is a class function on $G$ with values in $\ZZ$ and does not depend on
the choice of $x$.  By \cite{Huber01}, Theorem 4.1, $\sw_y$ is a virtual
character of $G$. Therefore, 
\[
    \sw_y(\F) := \gen{\sw_y,\tau}_G
\]
is an integer, which we call the {\em Swan conductor} of $\F$ at $x$.
By \cite{Huber01}, Proposition 4.2 (ii), this definition depends only
on $\F$ but not on the chosen representation
$\F\cong(f_*\QQb_\ell)[\tau]$.

It will be useful to have a formula for $\sw_y(\F)$ in terms of the
{\em break decomposition} of the $G$-module $V$ induced by the
filtration $(G^\gamma)$ of the group $G_x$,
\[
      V = \bigoplus_{\gamma\in\Gamma_x\otimes\QQ} V(\gamma).
\]
Here $V(\gamma)\subset V$ is defined for $\gamma<1$ as the subset of
$\gen{v-\sigma(v)\mid v\in V,\,\sigma\in G^\gamma}$ consisting of elements which
are invariant under $\cup_{\delta<\gamma}G^\delta$. Furthermore, $V(1):=V^P$,
where $P=\cup_{\gamma<1}G^\gamma$ is the maximal $p$-subgroup of $G_x$.  An
element $\gamma\in\Gamma_x\otimes\QQ$ is called a {\em break} for $\rho$ if
$V(\gamma)\not=1$. It is clear that a break is either a jump or equal to $1$. By
\cite{Huber01}, Corollary 8.4 we have
\begin{equation} \label{swanbreakeq}
  \sw_y(\F) \,=\, \sum_{i=1}^l \#\gamma_i\cdot\dim V(\gamma_i).
\end{equation}

The definition of $\delta(\F)$ is analogous. Let $\log_q:\RR_{>0}\to\RR$
be the logarithm to the basis $q:=\abs{\varpi^{-1}}$, where $\varpi$
is a uniformizer of the discretely valued field $K_0$. We set
\[
   \delta_x(\sigma) := -\abs{G_x}\log_qh^\flat_x(\sigma)\;\;\;
         \text{if $\sigma\neq 1$}, \quad
   \delta_x(1) := -\sum_{\sigma\neq 1} \delta(\sigma)
\]
and
\[
    \delta_y:= \Ind_{G_x}^G \delta_x.
\]
Finally, we define the discriminant conductor of $\F$ at $y$ as follows:
\[
      \delta_y(\F) := \gen{\delta_y,\tau}_G.
\]
Similar to \eqref{swanbreakeq}, we have a formula which
computes $\delta_y(\F)$ in terms of the break decomposition of
the $G$-module $V$:
\begin{equation} \label{swanbreakeq2}
  \delta_y(\F) = \sum_{i=1}^l -\log_q\gamma_i^\flat\cdot\dim V(\gamma_i).
\end{equation}
In particular, $\delta_y(\F)$ is a nonnegative rational number.

\subsection{Comparison with Huber's and Ramero's theory}

Let $\Xs$ be an open analytic curve with field of definition $K$. A
{\em compactification} of $\Xs$ is an open embedding
$\Xs\subset\Xs_1$, where $\Xs_1$ is an open analytic curve over $K$
such that $\Xs$ is contained in an affinoid subdomain of $\Xs_1$. Let
$\Xs_1\ad$ denote the analytic adic space associated to $\Xs_1$, see
\cite{Huber96}. We will associate to every end $x\in\partial\Xs$ a
certain point $x\ad\in\Xs_1\ad$.

Let $\As\subset\Xs$ be an open annulus representing the end $x$. After
enlarging the field $K$, if necessary, there exists a formal model
$\X$ of $\Xs_1$ with reduced special fiber $\X_s$, a closed subset
$Z\subset\X_s$ and a closed point $z\in Z$ such that $\Xs=]Z[_{\X}$
and $\As=]z[_{\X}$. Since $\As$ is an open annulus, $z$ is an ordinary
double point of $\X_s$. One of the irreducible components of $\X_s$
passing through $z$, say $W$, is not contained in $Z$; otherwise,
$\As$ would be contained in an affinoid which is itself contained in
$\Xs$, contradicting Definition \ref{opendef2}. Let $\Spm
A_K\subset\Xs_1$ be an affinoid subdomain containing $\Xs$. The
component $W$ gives rise to a discrete valuation
$\abs{\,\cdot\,}_\eta\ad$ on the fraction field of the affinoid
algebra $A_K$. The residue field of this valuation can be identified
with $k(W)$, the function field of $W$. By the definition of adic
spaces, $\abs{\,\cdot\,}_\eta\ad$ corresponds to a point
$\eta\in\Xs_1\ad$. In terms of the classification of points of
$\Xs_1\ad$ in \cite{Huber01}, \S 5, $\eta$ is a point of {\em type
II}. Such a point has infinitely many proper specialization to points
of {\em type III}. In particular, let $\abs{\,\cdot\,}_x\ad$ denote
the rank two valuation on the fraction field of $A_K$ which is the
composition of the valuation $\abs{\,\cdot\,}_\eta\ad$ with the
valuation on $k(W)$ corresponding to the point $z\in W$. Then
$\abs{\,\cdot\,}_x\ad$ corresponds to a point $x\ad\in\Xs_1\ad$ which
is of type III. By the definition of $\Xs_1\ad$ as a topological
space, $x\ad$ is contained in the closure of the point $\eta$. Extend
the valuation $\abs{\,\cdot\,}_\eta\ad$ to $\Frac(A_K\otimes_K K\ac)$
in the obvious way, and let $K_x\ad$ denote the henselization.

\begin{prop} \label{Kxprop}
  There is a canonical injection $K_x\ad\inj K_x$ of valued fields
  which induces an isomorphism on the value groups of the valuations.
\end{prop}

\proof Let $\U=\Spf B$ be an affine open formal subscheme of $\X$
containing the point $z$. The inclusion $\U\otimes K\inj\Xs_1$
corresponds to a morphism $A_K\to B\otimes K$. On the other
hand, the formal completion of $\U$ in $z$ can be identified with the
minimal formal model of the annulus $\As$. In particular,
$\OO_{\As}^\circ=\OOh_{\U,z}$. Hence we obtain a field extension
\[
      \Frac(A_K\otimes_K K\ac) \to \Frac(\OO_{\As}^\circ\otimes_R K\ac).
\]
By construction this is an extension of valued fields inducing an
isomorphism of the valuation groups, with respect to the valuations
$\abs{\,\cdot\,}_x\ad$ on the left and the valuation
$\abs{\,\cdot\,}_x$ on the right. The field $K_x\ad$ is
defined as the henselization of field on the left, whereas the field $K_x$ is the
completion of the field on the right. Since $K_x$ is henselian, we
obtain the desired injection $K_x\ad\inj K_x$. 
\Endproof

Let $\F$ be an admissible sheaf on $\Xs$ and suppose that $\F$ extends
to an admissible sheaf $\F_1$ on $\Xs_1$. Let $x\in\partial\Xs$ be an
end and $x\ad\in\Xs_1\ad$ the corresponding adic point on $\Xs_1$. Let
$\eta\in\Xs_1\ad$ be the unique generalization of $x\ad$. According to
\cite{Huber01} and \cite{Ramero}, we can define a Swan conductor and a
discriminant conductor
\[
     \sw_{x\ad}(\F_1) \in \ZZ, \qquad \delta_\eta(\F_1)\in\RR.
\]
These are defined in the same manner as $\sw_x(\F)$ and
$\delta_x(\F)$, with the valued field $K_x$ replaced by $K_x\ad$.
Therefore, Proposition \ref{Kxprop} shows:

\begin{cor} \label{swancor}
  For every compactification $\Xs_1\supset\Xs$ and every extension
  $\F_1$ of $\F$ to $\Xs_1$ we have
  \[
         \sw_x(\F) = \sw_{x\ad}(\F_1), \qquad \delta_x(\F) = \delta_\eta(\F_1).
  \]
\end{cor}

\begin{rem}
  It is plausible that for every admissible sheaf $\F$ on $\Xs$ there
  exists a compactification $\Xs_1\supset\Xs$ and an extension $\F_1$
  of $\F$ to $\Xs_1$, but the author does not know how to prove this.
\end{rem}

\subsection{Continuity}

Fix an element $R\in\Gamma\cup\{0\}$, $R<1$. If $R\neq 0$ we set
\[
        \Xs := \Cs(R,1) = \{\,t \,\mid\, R<\abs{t}<1 \,\},
\]
which is an open annulus; if $R=0$ we let 
\[
       \Xs:=\Ds(0,1) = \{\,t \,\mid\, \abs{t}<1 \,\}
\]
be the standard open disk. For every $r\in\Gamma$ with $R<r\leq 1$ we
set $s:=-\log_q r$ and define an open subset
\[
       \Xs_s := \{\,t\in\Xs_s \,\mid\, \abs{t}<r \,\},
\]
which is again an open annulus or an open disk. We let
$x_s\in\partial\Xs_s$ denote the `exterior' end corresponding to the
family of annuli $C(r',r)$ for $R<r'<r$. If $R\neq 0$ we shall identify the
`interior' end of $\Xs_s$ for all $s$ and denote it by $y$.

Let $\F$ be an admissible sheaf on $\Xs$. For $s=-\log_q r$ as above we define
\[
     \delta_\F(s) := \delta_{x_s}(\F|_{\Xs_s}),\qquad 
     \sw_\F(s) := \sw_{x_s}(\F|_{\Xs_s}).
\]

\begin{prop} \label{deltaprop}
  The association $s\mapsto \delta_\F(s)$ extends to a continuous and
  piecewise linear function $\delta_\F:[0,-\log_qR]\to\RR_{\geq 0}$.
  Similarly, the association $s\mapsto \sw_\F(s)$ extends to a right
  continuous and piecewise constant function
  $\sw_\F:[0,-\log_qR]\to\ZZ$. Furthermore:
  \begin{enumerate}
  \item
    The function $\delta_\F$ is convex.
  \item
    The function $\sw_\F$ is decreasing.
  \item
    For all $s\in[0,-\log_qR)$ we have
    \[
           \frac{\partial}{\partial s}\delta_\F(s+) = -\sw_\F(s).
    \]
  \item
    If $R=0$ then $\delta$ is decreasing and eventually zero. 
  \end{enumerate}
\end{prop}

\proof If $s>0$ then $\Xs$ is a compactification of $\Xs_s$ and hence
the end $x_s\in\partial\Xs_s$ corresponds to an adic point
$x\ad\in\Xs\ad$. Using Corollary \ref{swancor}, one can therefore
deduce Proposition \ref{deltaprop} for $s>0$ from results of
\cite{Ramero}, in particular Theorem 2.3.35 and Proposition 3.3.26.
At $s=0$ we do not have a definition of $\delta_\F$ and $\sw_\F$ in
terms of the adic space $\Xs\ad$, and we cannot directly use the
results of \cite{Ramero}. But all that remains to be shown is that
$\delta_\F$ and $\sw_\F$ are right continuous at $s=0$.

Let $f:\Ys\to\Xs$ be an \'etale $G$-Galois cover such that $f^*\F$ is
constant; then $\F\cong(f_*\QQb_\ell)[\tau]$ for some representation
$\tau$ of $G$. For $s\in(0,-\log_qR)$ and $r:=q^{-s}$ set $A_s:=\{t\in
\Xs \mid \abs{t}<r\}$ and let $\Bs\subset\Ys$ be a connected component
of $f^{-1}(\As)$. Using the semistable reduction theorem, one easily
shows that $\Bs$ is an open annulus for all $s$ sufficiently close to
zero (see \cite{open}). Therefore, the following lemma
implies that the function $\delta_\F$ (resp.\ $\sw_\F$) is linear
(resp.\ constant) on the interval $[0,s]$. The proof of the
proposition is now complete.  \Endproof

\begin{lem} \label{deltalem}
  Let $\Xs=\Cs(R,1)$ be an open annulus and $\Xs'=\Cs(R,r)$ a
  subannulus, with $R<r<1$. Let $x\in\partial\Xs$ be the `exterior'
  end corresponding to the family of annuli $C(r',1)$ for $R<r'<1$.
  Likewise, let $x'\in\partial\Xs'$ be the `exterior' end of $\Xs'$.
  Let $G$ be a finite group acting faithfully on $\Xs$ and fixing the
  end $x$. Then the following holds.
  \begin{enumerate}
  \item
    The action of $G$ fixes the annulus $\Xs'$ and the end $x'$.
  \item
    For all $\sigma\in G$, $\sigma\neq 1$ we have
    \[
        \#h_{x'}(\sigma) = \#h_x(\sigma)
    \]
    and 
    \[
        \log_q h_{x'}^\flat(\sigma) = 
          \log_q h_x^\flat(\sigma) + \#h_x(\sigma)\cdot\log_q r.
    \]
  \end{enumerate}
\end{lem}

\proof This follows from a direct computation, using Example
\ref{annulusexa} and \eqref{hsigmaeq}.  \Endproof

\begin{rem}
  The last part of its proof, including Lemma \ref{deltalem}, shows
  that most of Proposition \ref{deltaprop} (except maybe for the
  convexity of $\delta_\F$) is a rather straightforward consequence of
  the Semistable Reduction Theorem. This argument essentially goes
  back to L\"utkebohmert's paper \cite{Luetkebohmert93} on the
  non-archimedian Riemann Existence Theorem. See also
  \cite{Schmechta01}.
\end{rem}


\section{The boundary of Lubin-Tate space}

In this section we compute the filtration by higher ramification
groups at the boundary of the \'etale cover
$f_n:\Xs(\varpi^n)\to\Xs(1)$ of Lubin-Tate spaces of dimension one.

\subsection{Notation}

Let $F$ be a non-archimedian local field, e.g.\ a field which is
complete with respect to a discrete valuation $\abs{\,\cdot\,}$ and
whose residue field is finite, say with $q=p^f$ elements. We let $\OO$
denote the ring of integers of $F$. Furthermore, we choose a
uniformizer $\varpi$ of $\OO$. We assume that $\abs{\varpi}=q^{-1}$.

Let $K_0:=\Fh\nr$ denote the completion of the maximal unramified
extension of $F$, $R_0$ the valuation ring of $K_0$ and $k$ the
residue field of $K_0$. We choose an algebraic closure $K\ac$ of $K_0$
and extend the valuation $\abs{\,\cdot\,}$ from $F$ to $K\ac$. The
valuation group of $\abs{\,\cdot\,}$ is denoted by $\Gamma$. We write
$K$ to denote a finite extension of $K_0$ contained in $K\ac$ and $R$
for the valuation ring of $K$. Note that all this notation is
consistent with the notation used in Section 1.

Let $\Sigma_0$ be the unique formal $\OO$-module of height two over
$k$. Let $\X(1)=\Spf A_0$ be the universal deformation space of
$\Sigma_0$, and let $\Sigma\univ/A_0$ denote the universal deformation
of $\Sigma_0$. For each integer $n\geq 0$ be denote by
$\X(\varpi^n)=\Spf A_n$ the universal deformation space parameterizing
deformations of $\Sigma_0$ with a Drinfeld level structure of level
$\varpi^n$. We denote by
\[
      \phi_n: (\OO\varpi^{-n}/\OO)^2 \to (\Sigma\univ\otimes_{A_0}A_n)[\varpi^n]
\]
the tautological level-$\varpi^n$-structure on
$\Sigma\univ\otimes_{A_0}A_n$. For each $n\geq 0$ we let 
\[
        \Xs(\varpi^n) := \X(\varpi^n)\otimes_{R_0}K
\]
be the generic fiber of the formal scheme $\X(\varpi^n)$. We regard
$\Xs(\varpi^n)$ as an open analytic curve; recall that this means
essentially that the field of definition $K$ is a `sufficiently large'
finite extension of $K_0$. We remark that the right choice of $K$ will
in general depend on $n$. For instance, we want that the connected
components of $\Xs(\varpi^n)$ stay connected over any extension of
$K$. For this it is necessary and sufficient that the field $K$
contains the abelian extension of $K_0$ corresponding to the group of
$n$-units $1+\OO\varpi^n\subset\OO^\times$ via local class field
theory.

Write $G_n$ for the finite group $\GL_2(\OO/\varpi^n)$.  It is well
known that the natural map
\[
     f_n:\Xs(\varpi^n)\to\Xs(1)
\]
is an \'etale $G_n$-torsor. An element $\sigma\in G_n$ acts on
$\Xs(\varpi^n)$ by composing the tautological
level-$\varpi^n$-structure $\phi_n$ with $\sigma$, considered as a
linear automorphism $\sigma$ of $(\OO\varpi^{-n}/\OO)^2$.

\subsection{The ramification filtration}

By the fundamental result of Lubin-Tate and Drinfeld, the universal
deformation ring of $\Sigma_0$ is a power series ring over $R_0$,
i.e.\ $A_0=R_0[[T]]$. This means that we can identify the Lubin-Tate
space $\Xs(1)$ with the standard open unit disk.

Let $x\in\partial\Xs(1)$ be the unique end of $\Xs(1)$. Fix an integer
$n\geq 1$ and choose an end $y\in f_n^{-1}(x)=\partial\Xs(\varpi^n)$.
Let $K_y/K_x$ denote the corresponding extension of valued fields and
$G_y\subset G_n$ its Galois group. Let $h_0<\cdots<h_l$ denote the
jumps in the ramification filtration $(G_h)_h$ of the group $G_y$.
(See Section 1 for the notation.) As a convenient notational device,
we define the group homomorphism
\[
    -\Log:\Gamma_y\to\QQ\times\ZZ, \quad 
       (r,\,\gamma_y^i) \mapsto (-\log_q r, \,i).
\]

\begin{prop} \label{Gyprop}
  We may choose $y\in f_n^{-1}(x)$ such that 
  \[
      G_y \,=\, \Bigl\{\,\begin{pmatrix} 
                a & b \\ 0 & \;a^{-1} \end{pmatrix}
         \,\mid\, a\in\OO_{F,n}^\times,\, b\in\OO_{F,n}\,\Bigr\}.
  \]
  Then the jumps and the corresponding higher ramification groups are
  as follows.
  \begin{enumerate}
  \item
    For $i=1,\ldots,n-1$ we have
    \[
        -\Log h_i = (0,\,q^{2i}-1)
    \]
    and
    \[
        G_{h_i} \,=\, \Bigl\{\,\begin{pmatrix} 
            a & b \\ 0 & a^{-1} \end{pmatrix} 
            \,\mid\, a \equiv 1 \pmod{\varpi^i}\,\Bigr\}.
    \]
  \item
    For $j=0,\ldots,n-1$ and $i=n+j$ we have
    \[
        -\Log h_i = \bigl(\frac{1}{q^{n-j-1}(q-1)},\, -q^{2n-1}-1\bigr)
    \]
    and
    \[
        G_{h_i} \,=\, \Bigl\{\,\begin{pmatrix} 
                 1 & b \\ 0 & 1 \end{pmatrix} 
            \,\mid\, b \equiv 0 \pmod{\varpi^j}\,\Bigr\}.
    \]
  \end{enumerate}
\end{prop}

A proof of this proposition will be given in \S \ref{Gyproof} below. 

\begin{cor} \label{Gycor}
  For the jumps $\gamma_i=(\gamma_i^\flat,\#\gamma_i)$ in the upper
  numbering we have the following formulas:
  \begin{equation} \label{gammaieq2}\renewcommand{\arraystretch}{2}
    \#\gamma_i =\; \left\{\begin{array}{cl}
      (q+1)(1+q+\ldots+q^{i-1}), & i=1,\ldots,n-1 \\
      \displaystyle-\frac{q+1}{q-1},         & i=n,\ldots,2n-1,
                   \end{array}\right.
  \end{equation}
  and
  \begin{equation} \label{gammaflateq}
    -\log_q\gamma_{2n-1}^\flat = \frac{nq-n+1}{q-1}.
  \end{equation}
\end{cor}

\proof This follows from Proposition \ref{Gyprop} by a direct
computation, using \eqref{gammaieq}.  \Endproof

\subsection{The proof of Proposition \ref{Gyprop}} \label{Gyproof}

In this section we write $\Sigma=\Sigma\univ$ for the universal
deformation of the formal $\OO$-module $\Sigma_0$ over the ring
$A_0=R_0[[T]]$. We fix a parameter $X$ for the formal group law
underlying $\Sigma$. We write $[a](X)=aX+\ldots\in A_0[[X]]$ for the
endomorphisms of $\Sigma$ corresponding to an element $a\in\OO$. It is
well known that for a suitable choice of the parameters $X$ and $T$ we
have
\begin{equation} \label{picongruence}
   [\varpi](X) \,\equiv\, \varpi X+TX^q+X^{q^2} \pmod{\varpi X^q}.
\end{equation}
See e.g.\ \cite{HazewinkelFG}. For $n\geq 1$, let
\[
      \phi_n:\OO_n^2 \To \Sigma[\varpi^n]
\]
denote the tautological Drinfeld level-$\varpi^n$-structure over
$\X(\varpi^n)=\Spf A_n$. We identify the $\OO$-module $\Sigma[\varpi^n]$
with the subset of the maximal ideal of $A_n$ whose elements satisfy the
equation
\begin{equation} \label{pineq}
   [\varpi^n](X) = 
   \underbrace{[\varpi]\circ\cdots\circ[\varpi]}_{\text{$n$ times}}(X) \,=\, 0.
\end{equation}
The $\OO$-module structure on this set is given by the formal group law
$X+_{\Sigma}Y=X+Y+\ldots\in A[[X,Y]]$ and the formal power series
$[a](X)=aX+\ldots$. It is shown in \cite{DrinfeldEM} that
\begin{equation} \label{Aneq}
  A_n = A[u_n,v_n],
\end{equation}
where
\[
       u_n := \phi_n(1,0), \quad v_n:=\phi_n(0,1)
\]
is the standard $\OO_n$-basis of $\Sigma[\varpi^n]$ determined by
$\phi_n$. An element $\bigl(\begin{smallmatrix} a & b \\ c &
d\end{smallmatrix}\bigr)\in\GL_2(\OO_n)$ acts on $A_n$ over $A$ by the formula
\begin{equation} \label{abcdeq}
\begin{split}
   u_n &\longmapsto [a](u_n)+_{\Sigma}[c](v_n), \\
   v_n &\longmapsto [b](u_n)+_{\Sigma}[d](v_n).
\end{split}
\end{equation}

Let $H_n\subset\Sigma[\varpi^n]$ be the sub-$\OO$-module generated by
$u_n$.  Let $\Sigma_n':=\Sigma/H_n$ be the quotient (in the category of
formal $\OO$-modules over $A_n$) of $\Sigma$ under $H_n$. With respect to a
suitable parameter $X'$, the canonical morphism
$\Sigma\to\Sigma_n'$ is given by $X'=\alpha_n(X)\in A_n[[X]]$,
where
\begin{equation} \label{Serreseq}
     \alpha_n(X) \,=\, \prod_{a\in\OO_n} (X+_{\Sigma}[a](u_n)),
\end{equation}
see e.g.\ \cite{HazewinkelFG}. Let $\beta_n:\Sigma_n'\to\Sigma$ be the
`dual' isogeny to $\alpha_n$, i.e.\ the formal power series $\beta_n\in
A_n[[X']]$ such that
\begin{equation} \label{alphabetaeq}
     [\varpi^n](X) \,=\, \beta_n\circ\alpha_n(X).
\end{equation}

Let $x\in\partial\Xs(1)$ be the unique end of the disk $\Xs(1)$. An
end $y\in f_n^{-1}(x)$ corresponds to a valuation $\abs{\,\cdot\,}_y$
on $\Frac(A_n)$ extending $|\cdot|_x$. It follows from \eqref{Aneq}
that the field $K_y$ is the finite extension of $K_x$ generated by
$u_n,v_n$,
\begin{equation} \label{KyKxeq}
         K_y = K_x(u_n,v_n).
\end{equation}
We write $\val_y:=-\log_q\abs{\,\cdot\,}_y^\flat:K_y^\times\to\QQ$ for the
exponential version of the rank one valuation $\abs{\,\cdot\,}_y^\flat$. 

\begin{lem} \label{cansubgrouplem}
  We can choose $y\in f_n^{-1}(x)$ such that
  \[
          \val_y(u_n) \,=\, \frac{1}{(q-1)q^{n-1}}, \quad 
          \val_y(v_n) \,=\, 0.
  \] 
  If this is the case then $G_y$ is contained in the subgroup of
  $\SL_2(\OO_n)$ consisting of upper triangular matrices.
\end{lem}

\proof (Compare with \cite{Lubin67}.) We prove the lemma by induction on
$n$. Suppose first that $n=1$.  By \eqref{picongruence}, the Newton polygon of
$[\varpi](X)$ with respect to $\val_x$ has a unique negative slope $-1/(q-1)$
over the interval $[1,\ldots,q]$. It follows that $\Sigma[\varpi]$ has a
`canonical subgroup' $H$ whose nonzero elements $w$ satisfy
$\val_y(w)=1/(q-1)$. For $w\in\Sigma[\varpi]-H$ we have $\val_y(w)=0$. By
classical valuation theory, we can choose $y=y_1\in\Psi_1^{-1}(x)$ such that
$H=H_1=\gen{u_n}$. Then $G_y$ is contained in the stabilizer of $H_1$ which
consists of upper triangular matrices (use \eqref{abcdeq}). This proves the
case $n=1$.  The induction step from $n$ to $n+1$ is similar. For instance,
one uses the fact that $u_{n+1}$ is a solution of the equation
\[
     [\varpi](X)-u_n = 0.
\]
By \eqref{picongruence} and the induction hypothesis, the first slope of the
Newton polygon of this equation with respect to $\val_{y_n}$ is
$-1/q^n(q-1)$.  \Endproof

From now on we shall assume that $y$ is chosen as in Lemma
\ref{cansubgrouplem}. Set $\p_y$ denote the prime ideal of the
valuation ring $K_y^+$ of $K_y$ corresponding to the valuation
$\abs{\,\cdot\,}_y^\flat$. We write $K_y\sptilde:=\Frac(K_y^+/\p_y)$
for its residue field.  Lemma \ref{cansubgrouplem} and Equation
\eqref{Serreseq} show that we have the congruence
\begin{equation} \label{alphacongeq}
  \alpha_n(X) \,\equiv\, X^{q^n} \pmod{\p_y}.
\end{equation}
Using \eqref{picongruence}, \eqref{alphabetaeq} and induction on $n$ one
concludes that
\begin{equation} \label{betaEeq}
   \beta_n(X) \,\equiv\, E\circ E^{(1)}\circ\cdots\circ
                         E^{(n-1)}(X) \pmod{\p_y},
\end{equation}
where 
\[
     E^{(i)}(X) := T^{q^i}X+X^q.
\]
Set $w_n:=\alpha_n(v_n)$. Then $\beta_n(w_n)=0$ and by \eqref{alphacongeq} we
have $w_n\equiv v_n^{q^n}\pmod{\p_y}$. In particular, we have
$\val_y(w_n)=0$. Let $M_n:=K_x(w_n)$ denote the field extension of $K_x$
generated by $w_n$. We write $M_n\sptilde:=\Frac(M_n^+/(\p_y\cap M_n^+))$ for
the residue field of the valuation $\val_y|_{M_n}$. Note that the restriction of
$\#_y$ to $M_n$ induces a discrete valuation on $M_n\sptilde$. Let
$U_n\subset\SL_2(\OO_n)$ denote the subgroup of upper triangular, unipotent
matrices.

\begin{lem} \label{Gylem2}
\begin{enumerate}
\item
  The image of $w_n$ in the residue field 
  $M_n\sptilde$ is a uniformizer (with respect to the discrete valuation
  induced by $\#_y$). 
\item
  The field $M_n$ is the fixed field of $K_y$ of the subgroup $U_n\cap
  G_y\subset G_y$. 
\item
  The map 
  $\bigl(\begin{smallmatrix} a & b \\ 0 & a^{-1}\end{smallmatrix}\bigr)
  \mapsto a$ induces an isomorphism
  \[
      G_y/(G_y\cap U_n) \cong\OO_n^\times.
  \]
\end{enumerate}
\end{lem}

\proof
Let $z_n$ denote the image of $w_n$ in $M_n\sptilde$. By induction, we define
elements $z_m\in M_n\sptilde$ for $m=n-1,\ldots,0$ by
\[
   z_{m-1}:=E^{(m-1)}(z_m).
\]
By \eqref{betaEeq} we have $z_0\equiv\beta_n(w_n)\equiv 0\pmod{\p_y}$. Let
$L_m:=K_x\sptilde(z_m)\subset M_n\sptilde$ be the subextension of residue
fields generated by $z_m$. The element $z_1$ satisfies the equation over
$K_x\sptilde$ 
\[
    E(X)/X = X^{q-1} + T = 0.
\]
Since $T$ is a uniformizer of $K_x\sptilde$, we have $[L_1:K_x\sptilde]=q-1$
and $z_1$ is a uniformizer of $L_1$. For $m=1,\ldots,n-1$, the element
$z_{m+1}$ is a solution of the equation over $L_m$
\begin{equation} \label{zm1eq}
    E^{(m)}(X)-z_m= X^q + T^{q^m}X - z_m = 0.
\end{equation}
By induction, one proves that this is an Eisenstein equation and hence
irreducible over $L_m$ and that $z_{m+1}$ is a uniformizer of $L_{m+1}$. For
$m=n-1$, this gives Part (i) of the lemma. It also follows that
\begin{equation} \label{MnKxdegeq}
   [M_n:K_x]=[M_n\sptilde:K_x\sptilde]=(q-1)q^{n-1}=|\OO_n^\times|. 
\end{equation}
By the definition of $w_n$ and the isogeny $\alpha_n$, the
subfield $M_n\subset K_y$ is fixed under the action of $G_y\cap U_n$. On the
other hand, we have an injective homomorphism $G_y/(G_y\cap
U_n)\inj\OO_n^\times$. Now  \eqref{MnKxdegeq} implies Part (ii) and
(iii) of the lemma. 
\Endproof

\begin{lem} \label{Gylem3}
\begin{enumerate}
\item
  We have
  \[
      [K_y:M_n] = [K_y\sptilde:M_n\sptilde]=q^n=|\OO_n|.
  \]
\item
  The image of $v_n$ in $K_y\sptilde$ is a uniformizer (with respect to the
  discrete valuation induced by $\#_y$). 
\item
  The Galois group $G_y=\Gal(K_y/K_x)$ consists of all upper triangular
  matrices in $\SL_2(\OO_n)$.
\end{enumerate}
\end{lem}

\proof
By \eqref{alphacongeq} and the definition of $w_n$ we have
\[
         v_n^{q^n} \equiv w_n \pmod{\p_y}.
\]
Together with Lemma \ref{Gylem2} (i) it follows that (a) the field extension
$M_x(v_n)/M_x$ has degree $q^n$ (and its residue field extension
$M_x(v_n)\sptilde/M_n\sptilde$ is purely inseparable) and (b) $v_n$ is a
uniformizer of $M_n(v_n)\sptilde$. By Lemma \ref{Gylem2} (ii) we have
$[K_y:M_n]=|G_y\cap U_n]\leq |U_n|=q^n$. We conclude from (a) that
$K_y=M_n(v_n)$ and hence that Part (i) of the lemma holds. Also, Part (ii)
follows from (b) above. Finally, we have shown that $U_n\subset G_y$, and so
Lemma \ref{Gylem2} (iii) implies Part (iii) of the lemma.
\Endproof

\begin{lem} \label{Gylem4}
  For $m=1,\ldots,n$ we have
  \begin{align}
        \label{huneq}
     \val_y(u_m) & = \frac{1}{(q-1)q^{m-1}}, & \#_y u_m & = -q^{2n-1} \\
        \label{hvneq}
     \val_y(v_n) & = 0,                      & \#_y v_m & = q^{2(n-m)}.
  \end{align}
\end{lem}

\proof Let $y_m$ denote the image of $y$ in $\Psi_m^{-1}(x)$. Then
$\val_{y}|_{K_{y_m}}=\val_{y_m}$. Therefore, the formulae for $\val_y(u_m)$
and $\val_y(v_m)$ follow from Lemma \ref{cansubgrouplem} and the choice of
$y$. We also have
\begin{equation} \label{KyKymeq}
    \#_y|_{K_{y_m}} = [K_y:K_{y_m}]\cdot\#_{y_m} = 
       q^{2(m-n)}\cdot\#_{y_m},
\end{equation}
by Lemma \ref{Gylem3} (iii). Therefore, the formula for $\#_y(v_m)$ follows
from Lemma \ref{Gylem3} (ii). Moreover, in order to prove the formula for
$\#_y(u_m)$ it suffices to show that $\#_y(u_n)=-q^{2n-1}$. 

We proceed by induction on $n$. For $n=1$ we note that $u_1$ satisfies the
equation $[\varpi](X)=0$ over $K_x$. Substituting $X=\varpi^{1/(q-1)}Y$ and
using \eqref{picongruence}, we see that the image of $\varpi^{-1/(q-1)}u_1$ in
$K_x\sptilde$ is a solution to the equation
\[
    Y + TY^q = 0.
\]
We conclude that
\[
    q\cdot\#_{y_1}u_1 + \#_{y_1} T = \#_{y_1}u_1,
\]
which implies $\#_{y_1}u_1=-q$ and proves the claim for $n=1$. The induction
step from $n$ to $n+1$ is again similar: substitute $X=\varpi^{1/(q-1)q^{n-1}}Y$
into the equation $[\varpi](X)=u_n$ (of which $u_{n+1}$ is a solution).
\Endproof

We can now finish the proof of Proposition \ref{Gyprop}. If we choose $y$ as
in Lemma \ref{cansubgrouplem} then the first claim of the proposition (which
determines the group $G_y$) is proved by Lemma \ref{Gylem3} (iii). Let 
\[
    \sigma \;=\; \begin{pmatrix} a^{-1} & b \\ 0 & a \end{pmatrix}
\]
be an arbitrary element of $G_y-\{1\}$. To prove the proposition, it will
suffices to compute $h_y(\sigma)\in\Gamma_y$. If $a\neq 1\in\OO_n$ then we set
$i:=\val_{\varpi}(a-1)\in\{0,\ldots,n-1\}$. If $a=1\in\OO_n$ then
we set $j:=\val_{\varpi}(b)$ and $i:=n+j$. We claim that
\begin{equation} \label{hsigmaeq3}
    -\Log h_y(\sigma) \;=\; \left\{\begin{array}{cl}
      (0,\,0),  &  \text{for $i=0$}, \\
      (0,\,q^{2i}-1), & \text{for $i=1,\ldots,n-1$,} \\
      (\dfrac{1}{(q-1)q^{n-j-1}},\,-q^{2n-1}-1), & 
           \text{for $i\geq n$.}
                    \end{array}\right.
\end{equation}
Clearly, this claim implies the proposition.

By \eqref{hsigmaeq}, Lemma
\ref{Gylem3} (ii) and \eqref{abcdeq}, we have
\begin{equation} \label{hsigmaeq2}
  h_y(\sigma) \,=\, \Bigl\lvert\frac{\sigma(v_n)-v_n}{v_n}\Bigr\lvert_y
  \,=\, \Bigl\lvert\frac{([b](u_n)+_{\Sigma}[a](v_n))-v_n}{v_n}\Bigr\lvert_y.
\end{equation}
If $a\not\equiv 1\pmod{\varpi}$ then $[a](v_n)=av_n+(\text{higher
terms})\equiv av_n\not\equiv 0\pmod{\p_y}$ and $[b](u_n)\equiv 0\pmod{\p_y}$,
by \eqref{huneq} and \eqref{hvneq}. From \eqref{hsigmaeq2} we conclude that
\[
    -\Log h_y(\sigma) = -\Log \abs{a}_y = (0,\,0).
\]
This proves \eqref{hsigmaeq3} for $i=0$. Suppose now that
$i\in\{1,\ldots,n-1\}$ and write $a=1+\varpi^ic$. Then
\[
      [a](v_n) = [1+c\varpi^i](v_n) = v_n+_\Sigma[c](v_{n-i}) =
         v_n + cv_{n-i} + \ldots,
\]
where the remaining terms are units in $(K_y^+)_{\p_y}$ whose image in
$K_y\sptilde$ have valuation $>\#_yv_{n-i}=q^{2i}$, by \eqref{hvneq}. 
It follows that
\[
   -\Log h_y(\sigma) = -\Log\Abs{\frac{[a](v_n)-v_n}{v_n}}
                     = -\Log \Abs{\frac{cv_{n-i}}{v_n}}
                     = (0,\,q^{2i}-1).
\]
This proves \eqref{hsigmaeq3} for $i=1,\ldots,n-1$. Finally, for $i\geq n$ we
write $b=\varpi^jc$ and get
\[
    \sigma(v_n) = v_n +_\Sigma [c](u_{n-j}) = v_n + cu_{n-j} + \ldots,
\]
where the  remaining terms have $\p_y$-valuation
$>\val_y(u_{n-j})=1/(q-1)q^{n-j-1}$. Using \eqref{huneq} we conclude that
\[
  -\Log h_y(\sigma) = -\Log\Abs{\frac{cu_{n-j}}{v_n}} 
                    = \bigl(\frac{1}{(q-1)q^{n-j-1}},\, -q^{2n-1}-1\bigr).
\]
This completes the proof of Proposition \ref{Gyprop}.
\Endproof


\section{Ramification of supercuspidal 
          representations} \label{supercuspidal}

In this section we apply the results of the previous section to
compute the Swan and the discriminant conductor of the sheaves on the
Lubin-Tate space corresponding to the types of supercuspidal
representations of $\GL_2(F)$. We also draw some conclusions
concerning the cohomology of the Lubin-Tate tower.

\subsection{The cohomology of the Lubin-Tate tower} \label{supercuspidal1}

We continue with the notation introduced in the last section. We
also choose a prime number $\ell$ which is strictly bigger than
$p$. Then the order of the finite groups $G^{(n)}=\GL_2(\OO/\p^n)$ are
all prime to $\ell$. We write $K:=\GL_2(\OO)$ and $G:=\GL_2(F)$.  Let
$W_F$ denote the Weil group of $F$ and $I_F\subset W_F$ the inertia
subgroup.

Fix $n\geq 0$ and $i\in\{0,1,2\}$. We let
$H^i(\Xs(\varpi^n),\QQb_\ell)$ denote \'etale cohomology of the rigid
analytic space $\Xs(\varpi^n)\otimes_{K_0}\Kh\ac$, in the sense of
Berkovich \cite{Berkovich93}. See also \cite{open}, \S 2.1. Similarly,
one has cohomology with compact support,
$H^i_c(\Xs(\varpi^n),\QQb_\ell)$. We define
$H^1\para(\Xs(\varpi^n),\QQb_\ell)$, the {\em parabolic cohomology} of
$\Xs(\varpi^n)$, as the image of the natural map
$H^1_c(\Xs(\varpi^n),\QQb_\ell)\to H^1(\Xs(\varpi^n),\QQb_\ell)$. This
is a finite dimensional $\QQb_\ell$-vectorspace, together with a
continuous action of the group $G^{(n)}\times\OO_B^\times\times I_F$.

Set
\[
        \HH_0:= \varinjlim_n H^1\para(\Xs(\varpi^n),\QQb_\ell).
\]
This is an infinite-dimensional vector space with a continuous action
of the group $K\times\OO_B^\times\times I_F$. This action
extends, in a natural way, to an action of a certain subgroup
\[
     (G\times B^\times\times W_F)_0\subset G\times B^\times\times W_F.
\]
This subgroup is the kernel of the homomorphism $G\times
B^\times\times W_F\to\ZZ$ which sends $(g,b,\sigma)$ to the normalized
valuation of $\det(g)^{-1}N(b){\rm cl}(\sigma)$. Here $N:B^\times\to
F^\times$ is the reduced norm and ${\rm cl}:W_F\to F^\times$ the
inverse reciprocity map. See e.g.\ \cite{Fargues04} or \cite{LT}. We
let $\HH$ denote the representation of $G\times B^\times\times W_F$
induced from $\HH_0$.

\begin{thm}[Carayol]  \label{Carayolthm}
\begin{enumerate}
\item
  Let $\pi$ be an irreducible supercuspidal representation of $G$ over
  the field $\QQb_\ell$. As a representation of $B^\times\times W_F$ we have
  \[
       \Hom_G(\pi,\HH) \cong {\rm JL}(\pi)\spcheck\otimes{\rm L}(\pi)',
  \]
  where ${\rm JL}(\pi)$ is the image of $\pi$ under the local
  Jacquet-Langlands correspondence and ${\rm L}(\pi)$ is the image
  under the {\em Hecke}-correspondence (a certain normalization of the
  local Langlands correspondence). 
\item
  If $\pi$ is a smooth admissible irreducible representation of $G$
  which is not supercuspidal then $\Hom_G(\pi,\HH)=0$.
\end{enumerate}
\end{thm}

\proof Part (a) of this theorem is proved in \cite{Carayol86}, see
also \cite{Carayol90}. Part (b) is certainly known to the experts, but
there seems to be no explicit reference in the literature. We shall
deduce Part (b) from our computations of Swan conductors, see Corollary
\ref{sccor}.  \Endproof

Carayol has conjectured \cite{Carayol90} that Theorem \ref{Carayolthm}
extends to the group $\GL_n(F)$ and the corresponding Lubin-Tate
spaces of dimension $n-1$ for all $n\geq 2$. This conjecture has been
proved by Harris and Taylor \cite{HarrisTaylor}, along with the local
Langlands correspondence for $\GL_n$. Their method is a generalization
of Carayol's method \cite{Carayol86}, which in turn generalizes
arguments of Deligne \cite{DeligneLetter}. These arguments are quite
indirect. 

In a work in preparation, the author intends to give a new and more
direct proof of Theorem \ref{Carayolthm} which relies on an analysis
of the stable reduction of the spaces $\Xs(\varpi^n)$, studied in
\cite{LT}. The results of this paper, in particular Theorem \ref{thm1}
and Theorem \ref{thm2}, are a crucial ingredient for this analysis.

\subsection{The type of a supercuspidal representation}  \label{supercuspidal2}

In the following, all representations are defined over the
field $\QQb_\ell$. Recall that we are concerned with the groups
$G:=\GL_2(F)$, $K:=\GL_2(\OO)$ and $G^{(n)}=\GL_2(\OO/\p^n)$. We let
$K_n\subset K$ denote the principal congruence subgroup modulo
$\wp^n$. Let $U=U_0\subset K$ and $U_n\subset K_n$ be the subgroups
containing the upper triangular, unipotent matrices. We write
$H^{(n)}\subset G^{(n)}$ for the image of a subgroup $H\subset K$.
Let
\[
     K' \,:=\, \bigl\{\,\begin{pmatrix} a & b \\ c& d \end{pmatrix} \in K\,\mid\,
           c\equiv 0 \mod{\wp}\,\bigr\}
\]
denote the Iwahori subgroup of $K$ and set, for $n\geq 1$,
\[
     K_n' \,:=\,  1 + \begin{pmatrix} \wp^{n_2} & \wp^{n_1} \\ 
                                      \wp^{n_1+1} & \wp^{n_2} \end{pmatrix}
         \,\subset\, K,
\]
where $n_1:=[n/2]$ and $n_2:=[(n+1)/2]$ and where $\wp^0:=\OO$. Note that
$K_n'\subset K'$ is a normal subgroup for all $n$. We also let $Z$ (resp.\
$Z'$) denote the cyclic subgroup of $G$ generated by the element $\varpi\in
F^\times\subset G$ (resp.\ by the matrix $\Pi'$), where
\[
        \Pi' \,:=\, \begin{pmatrix} 0 & 1 \\ \varpi & 0 \end{pmatrix}.
\]
Note that $Z'$ normalizes $K'$ but not $K$. 

Let $\tau$ be a smooth and irreducible (and hence finite-dimensional)
representation of $K$. The {\em $K$-level} of $\tau$ is the minimal
integer $n\geq 1$ such that the restriction of $\tau$ to $K_n$ is
trivial. The {\em $K$-defect} of $\tau$ is the integer $n-r$, where
$r$ is the minimal integer such that the restriction of $\tau$ to
$U_r$ has a non-zero fixed vector. Using the filtration $(K_n')$, we
define in a similar way the $K'$-level and the $K'$-defect for a
smooth irreducible representation of $K'$.

A smooth and irreducible representation $\tau$ of $K$ is called {\em
  minimal} if its $K$-level cannot be lowered by twisting $\tau$ with
a one-dimensional character. The {\em minimal level} of $\tau$ is the
$K$-level of the twist of $\tau$ which is minimal.

Let $\pi$ be an irreducible admissible supercuspidal representation of
$G$.  For short we say that $\pi$ is a {\em supercuspidal}. Then there
exists a subgroup $J\subset G$, which contains and is compact modulo
the center of $G$, and a finite-dimensional representation $\sigma$ of
$J$ such that $\pi$ is the compactly induced representation
$\cInd_J^G(\sigma)$.  Furthermore, the pair $(J,\sigma)$ can be chosen
in a very specific way, as follows. We distinguish two cases. In the
first case, $J$ is the group $ZK$. Then the supercuspidal $\pi$ is
called {\em unramified}. The restriction of $\sigma$ to $K$ is called
the {\em $K$-type} of $\pi$ and is denoted by $\tau$. In the second case
$J=Z'K'$ and $\pi$ is called {\em ramified}. Then the $K$-type of
$\pi$ is defined as the induced representation
\[
        \tau:=\Ind_{K'}^K(\sigma|_{K'}).
\]
In both cases, the pair $(J,\sigma)$ is called a {\em type} for $\pi$. 

We have the following fundamental result of
Kutzko \cite{Kutzko78}. 

\begin{prop} \label{typeprop}
  Suppose that the $K$-type $\tau$ of $\pi$ is minimal of level $n$.
  Then the following holds.
  \begin{enumerate}
  \item If $\pi$ is unramified then $\tau$ has $K$-defect zero.
  \item If $\pi$ is ramified, then the restriction of $\sigma$ to $K'$
    has $K'$-level $2n-2$ and $K'$-defect zero. Furthermore, $n\geq 2$.
  \end{enumerate}
\end{prop}

\subsection{The Swan and the discriminant conductor of a type} 
\label{supercuspidal3}

Let $\pi$ be a supercuspidal representation with $K$-type $\tau$. Let
$m$ be the $K$-level of $\tau$. Considering $\tau$ as a representation
of the finite group $G^{(m)}=K/K_m$, we can define the admissible
sheaf on $\Xs(1)$
\[
   \F:=(f_{m*}\QQb_\ell)[\tau],
\]
where $f_m:\Xs(\varpi^m)\to\Xs(1)$ is the \'etale $K/K_m$-torsor
defined in the previous section. Let $x\in\partial\Xs(1)$ be the
unique end of the disk $\Xs(1)$. The following theorem, which is our
first main result, computes the Swan conductor and the discriminant
conductor of $\F$ at $x$ in the unramified case.

\begin{thm} \label{thm1} 
  Suppose that $\pi$ is unramified. Let $n$ be the minimal $K$-level of $\tau$. 
  Then we have
    \[
         \sw_x(\F) = -\frac{q+1}{q-1}\cdot\dim\tau = -(q+1)q^{n-1}
    \]
  and
    \[
         \delta_x(\F) = \frac{nq-n+1}{q-1}\cdot\dim\tau =
            (nq-n+1)q^{n-1}.
    \]
\end{thm}

\proof By Proposition \ref{Gyprop}, the stabilizers of the boundary
components of $\Xs(\varpi^m)$ are contained in the subgroup of $G_m$
of elements of determinant one. Therefore, if $\tau'$ is a twist of
$\tau$ by a one-dimensional character (which factors through
the determinant) and $\F'$ is the sheaf corresponding to $\tau'$
then $\sw_x(\F)=\sw_x(\F')$ and $\delta_x(\F)=\delta_x(\F')$. We may
therefore assume that $\tau$ is minimal of level $n$. For the rest of
the proof, we consider $\tau$ as a representation of $G^{(n)}$. Let $V$ be
the vector space underlying $\tau$. It follows from the construction
of $\tau$ in \cite{Kutzko78} that $\dim V=(q-1)q^{n-1}$.

Let $x_n\in\partial\Xs(\varpi^n)$ be an end whose stabilizer
$G_{x_n}^{(n)}\subset G^{(n)}$ consists of upper triangular matrices.
Let $\gamma_1>\ldots>\gamma_{2n-1}$ be the jumps for the filtration of
higher ramification groups $G^\gamma\subset G_{x_n}^{(n)}$. By
Proposition \ref{Gyprop} the subgroup $U_{n-1}^{(n)}\subset G^{(n)}$
is contained in $G^\gamma$ for all $\gamma\geq \gamma_{2n-1}$.  Since
$\tau$ has $K$-defect zero (Proposition \ref{typeprop}), the action of
$U_{n-1}^{(n)}$ on $V$ has no fixed vector. It follows that the break
decomposition of $V$ has a unique break at $\gamma_{2n-1}$, i.e.\
$V=V(\gamma_{2n-1})$. Using \eqref{swanbreakeq}, \eqref{swanbreakeq2}
and Corollary \ref{Gycor} we get
\[
     \sw_x(\F) = \#\gamma_{2n-1}\cdot\dim V = -(q+1)q^{n-1}
\]
and 
\[
    \delta_x(\F) = -\log_q\gamma_{2n-1}^\flat\cdot\dim V =
       (nq-n+1)q^{n-1}.
\]
\Endproof

\begin{cor} \label{thm1cor}
  With notation as in Theorem \ref{thm1}, we have
  \[
        \dim H^1\para(\Xs(1),\F) = 2q^{n-1}.
  \]
\end{cor}

\proof The Ogg-Shafarevich formula in \cite{Huber01} gives 
\[
    \sum_{i=0}^2 (-1)^i\dim H^i_c(\Xs(1),\F) = {\rm rank}\,\F + \sw_x(\F)
       = - 2q^{n-1}.
\]
Clearly $H^0_c(\Xs(1),\F)=0$. Moreover, $\dim H^2_c(\Xs(1),\F)=\dim
H^0(\Xs(1),\F)$ is equal to the dimension of the space of fixed
vectors of the representation $\tau$. Since $\tau$ is irreducible and
nontrivial, this number is also zero. We conclude that $\dim
H^1_c(\Xs(1),\F)=2q^{n-1}$. It remains to show that the map
$H^1_c(\Xs(1),\F)\to H^1(\Xs(1),\F)$ is an isomorphism. Indeed, the
dimension of the kernel and of the kokernel of this map equals the
intertwining number of $\tau$ with the trivial representation of the
stabilizer $G_y^{(n)}$ of an end $y\in\partial\Xs(\varpi^n)$. But this
number is zero because $G_y^{(n)}$ contains $U^{(n)}$ and $\tau$ has
$K$-defect zero.  \Endproof

Let us now assume that the supercuspidal $\pi$ is ramified. Let
$(J,\sigma)$ be the type of $\pi$, with $J=Z'K'$. We write
$\Xs_0(\varpi):=\Xs(\varpi)/K'$ for the \'etale cover of $\Xs(1)$
corresponding to the subgroup $K'\subset K$. Let $\F'$ denote the
admissible sheaf on $\Xs_0(\varpi)$ corresponding to the restriction
of $\sigma$ to $K'$. Since the $K$-type $\tau$ of $\pi$ is the induced
representation of the restriction of $\sigma$ to $K'$, the pushforward
of $\F'$ to $\Xs(1)$ can be identified with the sheaf $\F$.

Let $y\in\partial\Xs(\varpi)$ be an end of $\Xs(\varpi)$ such that the
stabilizer $G_y^{(1)}\subset G^{(1)}$ of $y$ is equal to the group of
upper triangular matrices of determinant one (Proposition
\ref{Gyprop}).  We see that $\Xs_0(\varpi)$ has exactly two ends,
corresponding to the double cosets $G_y\backslash K/K'$. (It is not
hard to show that $\Xs_0(\varpi)$ is an open annulus, see \cite{LT}.)
Let $y_1\in\partial\Xs_0(\varpi)$ be the image of $y$ and let $y_2$ be
the other end of $\Xs_0(\varpi)$. The next theorem computes the Swan
and the discriminant conductor of the sheaf $\F'$ at the two ends
$y_1$ and $y_2$.

\begin{thm} \label{thm2} Suppose that the supercuspidal $\pi$ is
  ramified and that its $K$-type has minimal $K$-level $n$. Let $\F'$
  be the admissible sheaf on $\Xs_0(\varpi)$ induced from the type
  $\sigma$ of $\pi$.  Let $y_1,y_2$ be the two ends of
  $\Xs_0(\varpi)$. Then
  \[
     \sw_{y_1}(\F') = \sw_{y_2}(\F') = -\frac{q+1}{q-1}\cdot\dim\sigma
                                    = -(q+1)q^{n-2}
  \]
  and 
  \[
     \delta_{y_1}(\F') = \delta_{y_2}(\F') = \frac{nq-q-n}{q-1}\dim\sigma
                                   = (nq-q-n)q^{n-2}.
  \]
\end{thm}

\proof It is shown in \cite{LT} that the pair $(\Xs_0(\varpi),\F')$
has an automorphism which switches the two ends $y_1$ and $y_2$. It
follows that the Swan and the discriminant conductor of $\F'$ at the
two ends are equal. Therefore it suffices to compute $\sw_{y_1}(\F')$
and $\delta_{y_1}(\F')$. Note also that the dimension of $\sigma$ is
equal to $(q-1)q^{n-2}$, by its construction in \cite{Kutzko78}. By
Proposition \ref{typeprop}, the restriction of $\sigma$ to $K'$ has
$K'$-level $2n-2$ and $K'$-defect zero.  Therefore, the restriction of
$\sigma$ to $U_r$ contains a fixed vector if and only if $r\geq n-1$.
Backed up by all these preliminary remarks, the proof proceeds exactly
as the proof of Theorem \ref{thm1}.  \Endproof

\begin{cor} \label{thm2cor}
  Let $\F$ be the sheaf on $\Xs(1)$ corresponding to the $K$-type
  $\tau$ of a ramified supercuspidal representation $\pi$. Let $n$
  denote the minimal $K$-level of $\tau$. Then
  \[
         \dim H^1\para(\Xs(1),\F) = 2(q+1)q^{n-2}.
  \]
\end{cor}

\proof
We have already remarked that the sheaf $\F$ is the pushforward of the
sheaf $\F'$ from Theorem \ref{thm2} via the map
$\Xs_0(\varpi)\to\Xs(1)$. Therefore, we have
\[
        H^i(\Xs(1),\F) = H^i(\Xs_0(\varpi),\F')
\]
for all $i$, and the same holds for cohomology with compact support
and parabolic cohomology. The Ogg--Shafarevich formula from
\cite{Huber01}, applied to the open annulus $\Xs_0(\varpi)$, gives
\[
   \sum_{i=0}^2 (-1)^i \dim H^i_c(\Xs_0(\varpi),\F') = 
    \sw_{y_1}(\F') + \sw_{y_2}(\F') = -2(q+1)q^{n-2},
\]
by Theorem \ref{thm2}. For the proof that $H^i_c(\Xs_0(\varpi),\F')=0$ for
$i=0,2$ and that $H^1_c(\Xs_0(\varpi),\F')\cong
H^1(\Xs_0(\varpi),\F')$ one proceeds as in the proof of Corollary
\ref{thm1cor}.  \Endproof

\begin{rem} \label{thm1rem}
  Let $\HH$ be the representation of $G\times B^\times\times
  W_F$ defined in \S \ref{supercuspidal1}. Corollary \ref{thm1cor} and
  Corollary \ref{thm2cor}
  imply that for an irreducible supercuspidal $\pi$ we have
  \[
      \dim \Hom_G(\pi,\HH)\; =\;\; \left\{
       \begin{array}{cl}
          4q^{n-1},\;\; & \text{if $\pi$ is unramified,} \\
          2(q+1)q^{n-2},\;\; & \text{if $\pi$ is ramified,}
       \end{array} \right.
  \]
  where $n$ is the minimal level of $\pi$. Indeed, using Frobenius
  reciprocity and the definition of the type of $\pi$ we
  see that $\Hom_G(\pi,\HH)$ can be identified with (a) the direct sum
  of two copies of $\Hom_K(\tau,\HH_0)=H^1\para(\Xs(1),\F)$ if $\pi$
  is unramified and with (b)
  $\Hom_{K'}(\sigma,\HH_0)=H^1\para(\Xs_0(\varpi),\F')$ if $\pi$ is
  ramified. From Theorem \ref{Carayolthm} we conclude that $\dim {\rm
  JL}(\pi)=2q^{n-1}$ in the unramified and $\dim {\rm
  JL}(\pi)=(q+1)q^{n-2}$ in the ramified case. These formulas can be
  verified by looking at the explicit construction of supercuspidal
  representations of $B^\times$ by type theory. 
\end{rem}

\subsection{Only supercuspidals occur in the parabolic cohomology}

Fix a character $\epsilon:\OO^\times\to \QQb_\ell^\times$ and let $n$
denote the {\em exponent} of $\epsilon$, i.e.\ the minimal positive
integer $n$ such that $\epsilon$ is trivial on $1+\p^n$. Let
$K_0(n)$ be the subgroup of $K$ consisting of matrices whose
lower left entry is divisible by $\varpi^n$. Then 
\[
    \begin{pmatrix} a & b \\ c & d \end{pmatrix} \;\mapsto\;
      \epsilon(a)
\]
defines a character $\tilde{\epsilon}:K_0(n)\to\QQb_\ell^\times$. 
Let 
\[
      u=u_n(\epsilon) := \Ind_{K_0(n)}^K\;\tilde{\epsilon}
\]
denote the induced representation. One can show that $u$ is
irreducible and has $K$-level $n$ (see e.g.\ \cite{Casselman73}).  Let
$\F$ denote the admissible sheaf on $\Xs(1)$ corresponding to the
representation $u$.

\begin{prop} \label{scprop}
  We have $H^1\para(\Xs(1),\F)=0$.
\end{prop}

\proof Let $x\in\partial\Xs(1)$ be the unique end. Choose
$y\in\partial\Xs(\varpi^n)$ as in Proposition \ref{Gyprop}. In the
following, we shall freely use the notation from the statement of
Proposition \ref{Gyprop} and Corollary \ref{Gycor}. Let $V$ be the
vector space underlying the representation $u$. Note that $\dim
V=[K:K_0(n)]=(q+1)q^{n-1}$. Consider the break decomposition
\[
     V = \bigoplus_{i=1}^{2n-1} V(\gamma_i)
\]
of $V$ induced from the filtration of higher ramification groups at
$y$. It follows from Lemma \ref{sclem} below that $V(\gamma_i)=0$ for
$i=1,\ldots,n-2$ and 
\[
   \dim V(\gamma_i) = \;\left\{
     \begin{array}{cl}
       2, & i=n-1, \\
       q^{j+1}-q^j, & i=n+j,\,j=0,\ldots,n-2,\\
       q^n-1,       & i=2n-1.
     \end{array} \right.
\]
Then by \eqref{swanbreakeq} and Corollary \ref{Gycor} the Swan
conductor of $\F$ at $x$ is
\[\begin{split}
  \sw_y(\F) & =\; \sum_{i=1}^{2n-1} \#\gamma_i\cdot\dim V(\gamma_i) \\
            & =\; 2\,(q+1)(1+q+\ldots+q^{n-2}) 
                \,-\, \sum_{j=0}^{n-2} q^j(q+1) \,-\, (q^n-1)\frac{q+1}{q-1} \\
            & =\; (q+1)\,\big((1+\ldots+q^{n-2})-(1+\ldots+q^{n-1})\big) \\
            & =\; -(q+1)\,q^{n-1}.
\end{split}\]
Therefore, the Ogg-Shafarevich formula gives
\[
   \dim H^1_c(\Xs(1),\F) = -\sw_x(\F) - \dim V = 0.
\]
\Endproof

\begin{lem} \label{sclem}
  For $i=1,\ldots,n-1$ we have $V^{G^{\gamma_i}} =0$.
  For $j=0,\ldots,n-1$ and $i=n+j$ we have
  \[
        \dim V^{G^{\gamma_i}} = 1+q^j.
  \]
\end{lem}

\proof
Left to the reader.
\Endproof

\begin{cor} \label{sccor}
  Let $\HH$ be the representation defined in \S
  \ref{supercuspidal1}. Let $\pi$ be an irreducible smooth admissible
  representation of $G$. If $\pi$ is not supercuspidal then
  \[
       \Hom_G(\pi,\HH) = 0.
  \]
\end{cor}

\proof Clearly, if $\pi$ occurs in the $G$-representation $\HH$ then
the restriction of $\pi$ to the subgroup $K$ occurs in the
$K$-representation $\HH_0=\varinjlim
H^1\para(\Xs(\varpi^n),\QQb_\ell)$. If $\pi$ is not supercuspidal then
by \cite[Appendix]{HenniartAppendix}, $\pi|_K$ either contains the trivial
representation or a representation $u=u_n(\epsilon)$ for some
character $\epsilon$ of exponent $n$. But the trivial representation
does not occur in $\HH_0$ because $H^1\para(\Xs(1),\QQb_\ell)=0$ and
$u$ does not occur in $\HH_0$ by Proposition \ref{scprop}. This proves
the claim.
\Endproof

\begin{rem}
  Laurent Fargues has explained to me a much better proof of the
  statement of Corollary \ref{sccor}, which also works for Lubin-Tate
  spaces of arbitrary dimension. 
\end{rem}

\end{document}